\documentclass[preprint, 12pt, titlepage,  draft]{elsarticle_w}
\usepackage{amssymb,latexsym,amsfonts,amsmath,amsthm}
\usepackage{verbatim}

\usepackage{color}

{\makeatletter \@addtoreset{equation}{section}}

\newcommand{\Ker}{\mathop{\mathrm{Ker}}\nolimits}
\newcommand{\IIm}{\mathop{\mathrm{Im}}\nolimits}

\newcommand{\beq}{\begin{equation}}
\newcommand{\eeq}{\end{equation}}

\newcommand{\bpm}{\begin{pmatrix}}
\newcommand{\epm}{\end{pmatrix}}

\def\bmat{\left( \begin{array}}
\def\emat{\end{array} \right)}

\newcommand{\Ra}{{\Rightarrow}}

\newcommand{\h}{ \mathop{ \mathrm{h} {} }\nolimits }
\newcommand{\e}{ \mathop{ \mathrm{e} {} }\nolimits }

\newcommand{\N}{\mathbb{N}}

\newcommand{\Aut}{\mathop{\mathrm{Aut}}\nolimits}
\newcommand{\End}{\mathop{\mathrm{End}}\nolimits}

\newtheorem{theorem}{Theorem}[section]
 \definecolor{msk}{rgb}{0.65,0.3,0.3}

\definecolor{msk2}{rgb}{0.75,0.3,0.1}
\definecolor{r}{rgb}{0.75,0.1,0.4}

\definecolor{exp2}{rgb}{0.85,0.6,0.1}

\definecolor{exp4}{rgb}{0.85,0.3,0.5}

\definecolor{vvs}{rgb}{0.5,0.4,0.7}

\definecolor{vvs}{rgb}{0.5,0.4,0.7}

\definecolor{kf}{rgb}{0.9,0.0,0.1}
\definecolor{gr}{rgb}{0.3,0.0,0.7}

\begin{document}

%%\setpagewiselinenumbers
%%\modulolinenumbers[5]
%% \linenumbers

\journal{Linear Algebra and its Applications} 

\begin{frontmatter}

\title{Hyperinvariant subspaces of locally nilpotent linear 
transformations}

\begin{abstract}

A subspace $X$ of a vector space over a field $K$ 
is hyperinvariant  with respect to an endomorphism $f$ of $V$ 
if it  is invariant for all endomorphisms of $V$ that
commute with $f$. 
We assume that $f$ is locally nilpotent, that is, 
 every $ x \in V $ is annihilated by some power of $f$, 
and that $V$ is an infinite direct sum  of $f$-cyclic subspaces. 
In this note we  describe the lattice of  hyperinvariant
subspaces of $V$. 
We  extend    results
of Fillmore, Herrero and Longstaff  (Linear Algebra Appl.~17 (1977), 125--132) 
to infinite dimensional spaces. 
\end{abstract}


\author[label1]{Pudji Astuti}
\ead{pudji@math.itb.ac.id} 

\author[label2]{\corref{crr}Harald K. Wimmer}
\ead{wimmer@mathematik.uni-wuerzburg.de} 
\cortext[crr]{Corresponding author}

\address[label1]{Faculty of Mathematics
 and Natural Sciences,
Institut Teknologi Bandung,
Bandung 40132,    Indonesia}

\address[label2]{                  
Mathematisches Institut,
Universit\"at W\"urzburg,
97074 W\"urzburg, Germany
}

\fntext[label1]
{The work of the first author was supported by the program  ``Riset
dan Inovasi KK ITB''  of  the  Institut Teknologi Bandung.
}

\end{frontmatter}

%%%%%%%%%%%% 
%%%%%
%%%%%
%%%%%

%%%   \newpage   \tableofcontents     \linenumbers

\section{Introduction} \label{sct.intrd}

Let $V$ be a vector space over a field $K$ and let $f$ be  an
endomorphism of~$V$. 
A subspace $X$ of $V$ is called {\em{hyperinvariant}} (with
respect of $f$)  if $X$  is invariant for all endomorphisms of $V$ that
commute with $f$ (see \cite[p.~227]{Bb}, \cite[p.\ 305]{GLR}).
 An endomorphism $f$ of $V$  is said to be {\em{locally nilpotent}} \cite[p.~37]{Kpl}
 if every $ x \in V $ is annihilated by some power of $f$.  
%%% ? Halmos 
In this  note  we are concerned with 
locally nilpotent  endomorphisms
with the property that the underlying vector space 
$V$ is an infinite  direct sum  of finite-dimensional  $f$-cyclic subspaces.
It is the purpose of our paper  to describe the lattice of  hyperinvariant
subspaces of $V$.  We  extend   results
of Fillmore, Herrero and Longstaff 
 \cite{FHL} (see \cite[Chapter~9]{GLR}) 
to infinite-dimensional spaces. 

Let  $t_j  $, % $t_j < t_{j+1} $, 
$  j   \in \N  $, 
 be the dimensions of the 
cyclic subspaces that are  direct summands of $V$.
Then  $V$ has a decomposition 
 \beq \label{eq.dmsr} 
V =  \bigoplus  _{j \in \N}   V_{t_j} , \:  \quad  % 
V_{t_j}  = 
\bigoplus  _{\sigma  \in S_j } V_{t_j  \sigma } 
\quad \text{where} \quad    V_{ t_j \sigma } \cong  K[s] / f ^{t_j} K[s]  .
\eeq  
In \eqref{eq.dmsr} the direct summands of dimension $t_j$ are gathered together in  
subspaces $V_{t_j} $, respectively. 
We  assume   
%%%% \beq \label{eq.ordd}   
%%%0 < t_1 < t_2 <  t_3 \cdots  %  t_j < t_{j+1 } < \cdots  
 % t_j < t_{j+1} ,  \,  j   \in \N ,  \eeq 
\beq \label{eq.ord}  
   t_j < t_{j+1} ,  \:   \,  j   \in \N .
\eeq   
The main result of the paper is the following. 
%% Theorem~\ref{thm.m} below.

\begin{theorem}  \label{thm.m} 
Suppose  $V$ is locally nilpotent with respect to $f$ and let 
\eqref{eq.dmsr} 
and \eqref{eq.ord} hold.  For a subspace  $ X  \subseteq  V$ 
the following statements are equivalent. 
\begin{itemize} 
\item[{\rm{(i)}}]   $X$ is hyperinvariant.
\item[{\rm{(ii)}}] %%The subspace 
$X$ is of the form 
\beq   \label{eq.gsm} 
X = \bigoplus  \nolimits  _{ j \in \N  } \,    f^{r_j } V_{t_j}    
\eeq 
with 
\beq     \label{eq.ins}  
0 \le r_j \le  t_j,  \: j \in \N,  \
\eeq 
and 
  \beq  \label{eq.ins3} 
 r_{j }  \:  \le \:  r_{ \ell} , 
\quad  
 t_j - r_{j} \:  \le \:  t_\ell-  r_{ \ell} , 
\quad  
if 
\quad     j  \le \ell.
\eeq 
\item[{\rm{(iii)}}] We have 
\beq \label{eq.imkr} 
X = \sum \nolimits  _{j \in \N} 
 \, 
%%f^{r_j} V
\big( \IIm f^{r_j} \,   \cap \,  \Ker f ^{t_j - r_j}  \big) 
%%V[f ^{t_j - r_j} ] 
\eeq 
with  $(r_j) _{j \in \N} $ satisfying
\eqref{eq.ins} and   \eqref{eq.ins3}.  
\end{itemize} 
\end{theorem}

The proof of the theorem will be given in Section~\ref{sct.oof}. 
We first introduce some notation and recall 
basic concepts and facts.  
Let $ x \in V $.   Define $ f^0 x  = x $.
The smallest nonnegative integer $\ell$
 with
$f^{\ell} x = 0$
is  called
the {\em{exponent}} of $x$. We write  $\e(x) = \ell$.
A  nonzero vector  $x $   is said to have \emph{height} $q $ if $x \in f^q V$
and $x \notin f^{q+1} V$.
In this case we write $\h(x) = q$. We set   $ \h ( 0 ) =   \infty $.
We remark  that a decomposition \eqref{eq.dmsr} 
implies that $V$  has no nonzero  elements of infinite height. 
In  particular, 
if $x \in \oplus _{j = 1  }  ^r  V_{t_j}$,  $x \ne 0 $, then
 \, $ \h(x)  <  t_r$.
%% corr  $ \h(x)  <  t_{j_r} $ 
Let $ Y $ be a subspace of $V$.  
We set 
 $  Y[ f^i ] =  \Ker \big(f^i _{|Y} \big) $. 
Thus, in the representation  \eqref{eq.dmsr} 
we have  %%% cases 
\beq 
\label{eq.slc} 
V_{t_j}[ f^{t_i}  ] =   V_{t_j}  \quad \text{if} \quad  j \le i ,
\quad 
\text{and} \quad   V_{t_j}[ f^{t_i}  ] = f^{t_j - t_i}  V_{t_j} 
\quad  \text{if} \quad  j > i .
\eeq 
We note that the decomposition \eqref{eq.dmsr} is unique up to
isomorphism.   
A corresponding result for abelian $p$-groups  in \cite[p.\ 89]{FuI} 
shows that the cardinality 
%$  \omega_j$ 
of $ S_j $ 
%%If $  \omega_j$ is the cardinality of $ S_j $ 
is given by   %\,$\omega_j = 
\,$ \dim \big( f^{t_j- 1} V /  f^{t_j } V \big)[f ]  $. 
Let
\[
\langle x  \rangle \,  = \,  {\rm{span}} \{ f^i x , \, i \ge 0 \}  \, 
=  \, 
\,  \{ a(f)  x ; \,  a(s) \in K[s] \}
\]
 be  the $f$-cyclic subspace  generated by $ x $.  
%% To $ B  \subseteq V$ we associate the subspace 
%% \marginpar{B?wg} 
%% \,$  \langle B \rangle  = 
%%  \sum _{b \, \in \,  B } \,  \langle \,  b  \, \rangle  $.
Adapting a definition of \cite{FuI} we call a vector 
$u$  a {\em{generator}} of $V$  
if $u\ne0$ and 
$  V =   \langle u\rangle \oplus V_2 $ for some invariant subspace 
$ V_2  $.
 %%\quad \text{for some}   \quad V_2 \in  {\rm{Inv}}( V,f) . 
%%
The set  $B_u = \{u, f u , \dots , f^{\e(u) -1} u \} $ is a Jordan basis
of  the direct summand 
$ \langle u \rangle  $ with respect to the restriction  $ f_{|\langle u\rangle }$. 
The invariant subspaces of $ \langle u\rangle  $ are 
$ \langle f^ q  u\rangle  $, $ q = 0, 1 , \dots \e(u) $.  %% \marginpar{stl} 
%%  $\langle x  \rangle $ has a direct $f$-invariant complement. 
%%It will be convenient to
Given the decomposition  \eqref{eq.dmsr} 
we  choose   a set of generators $\{u_{t_j  \sigma } , \, 
j \in \N ,   \sigma \in S_j \} $ %%  of $V$
%%\marginpar{stl} 
such that    
\beq \label{eq.gnu} 
 V = \bigoplus \nolimits _{j \in \N}  V_{t_j },  \,\,
V_{t_j }  = \bigoplus \nolimits _{\sigma \in S_j}
\langle  u_{t_j \sigma} \rangle
  \quad    \text{where} \quad   \e( u_{t_j  \sigma } ) = t_j .
 \eeq 
%%%%%%and the sequence of exponents   $ ( t_j) _{j \in \N} $ 
%%%% is ordered as in \eqref{eq.ord}.
Then $x \in V $  can be represented as a sum 
\beq \label{eq.rpln} 
x =
 \sum \nolimits  _{j\in \N,  \,\, \sigma   \in S_j  }   c_{ j \sigma }(f) \,   u_{ t_j  \sigma }, 
\; \;
 c_{ j \sigma }(s) \in K[s] ,  %%\,  c_{ j \sigma } \ne 0 , 
\eeq 
where only finitely many polynomials  $  c_{ j \sigma }(s) $ are different from
the zero polynomial. 
The representation \eqref{eq.rpln} 
is unique if \,$ \deg  c_{ j \sigma }(s)  < \e ( u_{ t_j  \sigma }  )
 = t_j$. 
Since 
$ \{   u_{ t_j \sigma }; j \in \N,  \sigma \in S_j \} $ is a basis of $V$ 
we have 
\beq \label{eq.mxpt} 
\e( x   ) = \max\{ \e\!\big(   c_{ j \sigma }(f)  u_{ t_j \sigma } \big) ; \; 
 c_{ j \sigma }(f)  u_{ t_j \sigma }    \ne 0 \} .
\eeq 
Let $ \pi _{t_j \sigma} $ %% : V \to  \langle  u_{t_j \sigma } \rangle    $
denote  the projection of $V$  on $ \langle  u_{t_j \sigma } \rangle   $ 
corresponding to \eqref{eq.gnu}. 
Thus, if $x$ is represented by \eqref{eq.rpln} 
 then 
\,$  \pi _{t_j \sigma} x =   c_{ j \sigma }(f)  \,   u_{ t_j  \sigma } $.
%%c_{ j \sigma }(s) \in K[s] .  
%
Let $\End(V) $ and $\Aut(V) $ be the endomorphism ring 
 and the automorphism group  of  $V$, respectively. 
If  $ \eta \in \End(V) $  then 
$  \e ( \eta x ) \le  \e (x) $. 
Suppose 
\beq \label{eq.nmjt}
     \eta  u_{ t_j \sigma }  
   =     \sum _{\ell\in \N , \, \tau  \in S_\ell  }  b_{ \ell \tau }(f)\,  u_{ t_\ell \tau  } ,
\, \, \, 
 b_{ \ell \tau }(s) \in K[s] .
\eeq 
Then \eqref{eq.mxpt} implies 
\,$
\e\!\big(   b_{ \ell \tau }(f)  u_{ t_\ell \tau  } \big) \: \le 
 \:  \e (  u_{ t_j \sigma }  ) \: 
  = \:    t_j$\,  if \,$b_{ \ell \tau }(f)  u_{ t_\ell \tau  }  \ne 0 $. 
Hence \,$f ^{t_j }  b_{ \ell \tau }(f)  u_{ t_\ell \tau  } = 0 $.
From $\e( u_{ \ell \tau } ) = t_\ell $ follows 
\,$ s^{t_\ell} \mid s^{t_j}  b_{ \ell \tau }(s)  $.
%%%Therefore we obtain a  condition for the coefficients 
%%$  b_{ \ell \tau }(s) $ in \eqref{eq.nmjt}.   Namely,
Hence,   if $ \ell >  j $ 
then 
 \,$ b_{ \ell \tau }(s) = s^{t_{\ell}    - t_j } a_{ \ell  \tau  }(s) $\, 
for some \,$ a_{ \ell  \tau  }(s)  \in K[s]$.  

\bigskip \bigskip \bigskip

\section{Proof of Theorem~\ref{thm.m}}
   \label{sct.oof} 

\begin{proof} 
(i) $\Ra$ (ii) \, 
Let $V$ be given by \eqref{eq.gnu}  and \eqref{eq.ord}. 
Suppose $X \subseteq V$ is hyperinvariant. 
Set $  X_{t_j \sigma } =   \pi _{t_ j \sigma } X $.   
Then 
\[
  X  \cap \langle  u_{t_j \sigma } \rangle  \, \subseteq  \,   X_{t_j \sigma }  \, 
\subseteq  \,    \langle   u_{t_j \sigma } \rangle  .
\]
Because of  $  \pi _{t_j \sigma }  \in   \End(V) $ we have 
 $  X_{t_j \sigma }  \subseteq X  $.  Hence
\,$
  X_{t_j \sigma }  = %%\pi _{t_j\sigma}  X = 
X \cap  \langle  u_{t_j \sigma } \rangle , 
$\, 
and therefore $  X_{t_j \sigma }  $
%% \subseteq \langle   u_{t_j \sigma } \rangle $ 
is an
invariant   subspace  contained in $  \langle   u_{t_j \sigma } \rangle $, 
and we conclude that
\[   
  X_{t_j \sigma }  = \langle   f^{r_{j \sigma } } u_{t_j \sigma } \rangle  
 \quad \text{for some} \: \:    r_{j \sigma} \:  \:    \text{with} \quad 
  0 \le  r_{j \sigma} \le t_j . 
\]
Then 
\[
X \: \subseteq  \: \bigoplus _{j \in \N, \, \sigma \in S_j}
 \big( \pi _{t_j \sigma}  X \big)  \: = \:  
\bigoplus _{j \in \N, \, \sigma \in S_j} 
\big(   X \cap   \langle
 u_{t_j\sigma }  \rangle  \big)   
 \: = \:  
\bigoplus _{j \in \N, \, \sigma \in S_j}    \langle
 f^{ r_{j \sigma} }  u_{t_j\sigma }   \rangle     \: \subseteq  \:  X
\]
implies
\beq \label{eq.wno} 
X = \bigoplus _{j \in \N} X_{t_j}   \quad \text{with} \quad 
X_{t_j} = X  \cap V_{t_j}  = \bigoplus _{\sigma \in S_j}  
%%\big(   G \cap  
  \langle f^{ r_{j \sigma} }    u_{t_j \sigma } \rangle , \, \,  
j\in \N.   
\eeq 
Let $ \alpha _{\sigma \rho} ^{\langle j \rangle } \in \Aut(A)  $ be the 
automorphism with
exchanges the generators  
 $  u_{t_j \sigma }   $ and $  u_{t_j \rho} $
in the sense that 
%% with the effect that?
\[
 \alpha _{\sigma \rho} ^{\langle j\rangle }   ( u_{t_j \sigma } ,  u_{t_j \rho} )
=   ( u_{t_j \rho},  u_{t_j \sigma } ) ,  \:\: {\rm{and}} \:\: 
 \alpha _{\sigma \rho} ^{\langle j\rangle }   ( u_{t_j \lambda }) = u_{t_j \lambda} 
\:\: {\rm{if}} \:\:  \lambda \notin \{ \rho, \sigma\} . 
\]
Then \,$ \alpha _{\sigma \rho} ^ {\langle j \rangle }   X =  X $\, and 
\,$\alpha _{\sigma \rho} ^{\langle j \rangle }  X_{t_j}  = X_{t_j} $. 
Therefore 
\eqref{eq.wno}  yields  $ r_{j \sigma} =  r_{j \rho} $  
for $\sigma, \rho \in S_j $. 
Hence 
\beq \label{eq.sgz} 
 X_{t_j} = f^{ r_j }  V_{t_j}  \quad \text{for some} \:\;    r_j 
\:\;
\text{with} \:\;   0 \le r_j\le t_j. 
\eeq

To show that  inequalities  \eqref{eq.ins3}  are satisfied 
we use suitable   endomorphisms 
$\eta _k^{\ell }  \in \End(V) $.  
  Let  
$ j, k, \ell \in \N$, $ \tau \in S_j $.
For each $\ell$ we fix an element 
$ \hat \mu _{\ell} $ of $S_{\ell}$. 
%%Defining
Let 
 $ \eta _k ^{\ell } $  be defined on the set of  generators  $\{ u_{t_j \tau} \}$  
of $V$  by  
\[
\begin{aligned} 
  \eta _{k } ^{\ell}  \,    u_{t_j \tau} &  = 0  \hspace{24mm}
 {\rm{if}} \:  \:    j \ne k ,  
\\
 \eta _{k } ^{\ell}  \,    u_{t_k\tau} & =
\begin{cases}  u_{t_\ell  \hat{\mu} _{\ell}    }   & {\rm{if}}  \:  \: k > \ell 
\\
  f^{t_\ell  -  t_k} \,  u_{t_ \ell \hat{\mu}_\ell  }   & {\rm{if}}  \: \: k < \ell .
\end{cases} 
\end{aligned} 
\]
Then
\[ 
\IIm \eta _{k } ^{\ell} \:   =  \:  \langle  u_{t_ \ell \hat{\mu}_\ell  } \rangle 
\big[ f^{t_k} \big]
\quad \text{and} \quad 
\Ker \eta _{k } ^{\ell}  \: =    \: \bigoplus \nolimits _{j \in  \N, \,  j \ne k } V_{t_j}. 
\]
Hence \eqref{eq.wno}  and  \eqref{eq.sgz}  
imply  
\[
\eta _{k } ^{\ell} X \:  = \:  \eta _{k } ^{\ell} X_{t_k}  \: 
=  \:   \eta _{k } ^{\ell} \,  f^{r_k} V_{t_k}  . 
\]
If $ k > \ell $  then
\[
\eta _{k } ^{\ell} \,  f^{r_k} V_{t_k}
 \:  = \: 
\langle  f ^{r_k}  u_{t_\ell  \hat{\mu}_\ell}   \rangle 
 \:
  \subseteq  
\: 
 V_{t_\ell} \cap X = 
X_{t_\ell} \:  =  \:  \:  f^{r_{\ell}}  V_{t_\ell}
\: =
\: 
\bigoplus 
\nolimits _{\tau \in S_\ell}  \langle  f ^{r_\ell} u_{\ell \tau } \rangle .
\]
Hence \,$ \langle  f ^{r_k}  u_{t_\ell  \hat{\mu}_\ell}   \rangle  \:  \subseteq  
\langle  f ^{r_\ell}  u_{t_\ell  \hat{\mu}_\ell}   \rangle $,
and we obtain \,$ r_k \ge r_{\ell} $.   If 
\,$ k < \ell  $\,  then
 \[
   \eta _{k } ^{\ell} f ^{r_k} V_{t_k}  \: 
= \: 
\langle    f^{r_k}  f^ {t_\ell - t_k}   u_{t_\ell  \hat{\mu}_\ell   } \rangle 
 \subseteq   X_{t_\ell}  \: = \: f^{r_{\ell}}   V_{t_\ell}
\]
implies  
$\langle    f^{ {r_k} + t_\ell - t_k} \,  u_{t_\ell  \hat{\mu}_\ell   } \rangle 
\subseteq \langle    f^{r_{\ell}} \,  u_{t_\ell  \hat{\mu}_\ell   } \rangle  $.
Thus we obtain
$r_k + t_\ell - t_k \ge r_{\ell} $, 
that is, 
$ t_\ell - r_{\ell} \ge t_k - r_k $. 

(ii) $\Ra$ (iii)  
Let   $X$  be a subspace of the form \eqref{eq.gsm}
with     \eqref{eq.ins}   and  \eqref{eq.ins3}. 
From   $ 0 \le r_j \le t_j$ follows 
\, $  f^{r_j } V_{t_j}  =
 V_{t_j} [ f^{t_j - r_j} ] \subseteq  V [ f^{t_j - r_j} ] $.
Therefore 
\,$  f^{r_j } V_{t_j}  \subseteq  f^{r_j } V \cap  V [ f^{t_j - r_j} ]  $,
which implies \, 
$ X \subseteq \sum _{j \in \N} \big(  f^{r_j } V \cap  V [ f^{t_j - r_j} ] \big)$. 

Now let $x \in f^{r_j } V \cap  V [ f^{t_j - r_j} ]$, $x \ne 0$. 
Then $ x =  f^{r_j } y $ for some nonzero $y \in V$ with $f ^{t_j} y = 0 $. 
Let $y$ be decomposed in accordance with  \eqref{eq.gsm}  such that
\[
 y =  
   y_1 + \cdots + y_n 
 \quad \text{with}
 \quad 
 y _{\nu} \in  V_{t_{j_{\nu}}}, \,   y_{\nu}  \ne 0 , \,  \nu = 1, \dots , n. 
\]
Then 
\,$
\max  \, % \{\e( y_{\nu} ) \} = \e(y) \le t_j . 
\{\e( y_{\nu} ) ;  \;  \nu = 1, \dots , n  \} = \e(y) \le t_j $, 
and therefore
\beq \label{eq.roa} 
\e( f^{r_j }  y_{\nu} )\le   t_j -  r_j . 
\eeq 
We have  $x \in X $ if \,$f^{r_j }  y _{\nu } \in X $, $  \nu = 1, \dots , n $. 
Let $ j_{\nu} \le j $. Then \,$ r_{ j_{\nu}} \le r_j $, and therefore 
\[
f ^{r_j} y_{\nu}  \in  f ^{r_j}  V_{t_{j_{\nu}}}  \subseteq 
    f ^{r_{j _\nu } }  V_{t_{j_{\nu}}}  \subseteq X . 
\]
Let    $ j  \le j_{\nu} $.  Then  \,$ t_j  - r_j  \le t_{ j_{\nu} }  - r_{ j_{\nu} } $
and  \eqref{eq.roa} imply
\[
 f^{r_j }  y_{\nu}  \in   V_{  t_{ j_\nu }}\big[  f^{  t_j -  r_j }\big]
\subseteq 
    V_{t_{ j_\nu }} \big[  f ^{  t_{j_\nu }   -   r_{j_\nu }   }\big]
 = f ^{  r_{j_\nu} }   V_{ t_{ j_\nu }} \subseteq X . 
\]
Thus we have shown that
\,$
x = \sum \nolimits _{\nu = 1} ^n f^{r_j} y_{\nu} \, \in X $,
which proves the inclusion
\[ 
f^{r_j } V \cap  V [ f^{t_j - r_j} ] \subseteq X  \quad \text{for all} \quad
j \in \N.
\]

(iii) $\Ra$ (i)  This is obvious, since  
$ \IIm f^p $ and  $\Ker f^q$  
 are hyperinvariant subspaces for all $ p, q \in \N_0$.
\end{proof}

\end{document}